\documentclass[12pt]{article}
\usepackage{e-jc}

\usepackage{amsmath, epsfig, cite}
\usepackage{amssymb}
\usepackage{amsfonts}
\usepackage{mathrsfs}
\usepackage{latexsym}
\usepackage{amsthm}

\numberwithin{equation}{section}

\title{A $q$-analogue of some binomial coefficient identities of Y. Sun}

\author{Victor J. W. Guo$^1$ and Dan-Mei Yang$^2$ \\
\small Department of Mathematics, East
China Normal University\\
\small  Shanghai 200062, People's Republic of China\\
\small \texttt{$^1$jwguo@math.ecnu.edu.cn,\quad  $^2$plain\_dan2004@126.com}
}

\date{\dateline{Dec 1, 2010}{Mar 24, 2011}}

\begin{document}
\maketitle

\begin{abstract}
We give a $q$-analogue of some binomial coefficient identities of Y. Sun
[Electron. J. Combin. 17 (2010), \#N20] as follows:
\begin{align*}
\sum_{k=0}^{\lfloor n/2\rfloor}{m+k\brack k}_{q^2}{m+1\brack n-2k}_{q} q^{n-2k\choose 2}
&={m+n\brack n}_{q}, \\[5pt]
\sum_{k=0}^{\lfloor n/4\rfloor}{m+k\brack k}_{q^4}{m+1\brack n-4k}_{q} q^{n-4k\choose 2}
&=\sum_{k=0}^{\lfloor n/2\rfloor}(-1)^k{m+k\brack k}_{q^2}{m+n-2k\brack n-2k}_{q},
\end{align*}
where ${n\brack k}_q$ stands for the $q$-binomial coefficient.
We provide two proofs, one of which is combinatorial via partitions.
\end{abstract}

\section{Introduction}
Using the Lagrange inversion formula, Mansour and Sun \cite{MS} obtained the following two binomial coefficient identities:
\begin{align}
\sum_{k=0}^{\lfloor n/2\rfloor}\frac{1}{2k+1}{3k\choose k}{n+k\choose 3k}
&=\frac{1}{n+1}{2n\choose n}, \label{eq:s1}\\
\sum_{k=0}^{\lfloor (n-1)/2\rfloor}\frac{1}{2k+1}{3k+1\choose k+1}{n+k\choose 3k+1}
&=\frac{1}{n+1}{2n\choose n}\quad (n\geqslant 1).  \label{eq:s2}
\end{align}
In the same way, Sun \cite{Sun} derived the following binomial coefficient identities
\begin{align}
\sum_{k=0}^{\lfloor n/2\rfloor}\frac{1}{3k+a}{3k+a\choose k}{n+a+k-1\choose n-2k}
&=\frac{1}{2n+a}{2n+a\choose n}, \label{eq:s3}\\
\sum_{k=0}^{\lfloor n/4\rfloor}\frac{1}{4k+1}{5k\choose k}{n+k\choose 5k}
&=\sum_{k=0}^{\lfloor n/2\rfloor}\frac{(-1)^k}{n+1}{n+k\choose k}{2n-2k\choose n}, \label{eq:s4}\\
\sum_{k=0}^{\lfloor n/4\rfloor}\frac{n+a+1}{4k+a+1}{5k+a\choose k}{n+a+k\choose 5k+a}
&=\sum_{k=0}^{\lfloor n/2\rfloor}(-1)^k{n+a+k\choose k}{2n+a-2k\choose n+a}. \label{eq:s5}
\end{align}

It is not hard to see that both \eqref{eq:s1} and \eqref{eq:s2} are special cases of \eqref{eq:s3}, and \eqref{eq:s4}
is the $a=0$ case of \eqref{eq:s5}. A bijective proof of \eqref{eq:s1} and \eqref{eq:s3}
using binary trees and colored ternary trees has been given by Sun \cite{Sun} himself. Using the same model,
Yan \cite{Yan} presented an involutive proof of \eqref{eq:s4} and \eqref{eq:s5}, answering a question of Sun.

Multiplying both sides of \eqref{eq:s3} by $n+a$ and letting $m=n+a-1$, we may write it as
\begin{align}
\sum_{k=0}^{\lfloor n/2\rfloor}{m+k\choose k}{m+1\choose n-2k}={m+n\choose n}, \label{eq:new1}
\end{align}
while letting $m=n+a$, we may write \eqref{eq:s5} as
\begin{align}
\sum_{k=0}^{\lfloor n/4\rfloor}{m+k\choose k}{m+1\choose n-4k}
=\sum_{k=0}^{\lfloor n/2\rfloor}(-1)^k{m+k\choose k}{m+n-2k\choose m}.  \label{eq:new2}
\end{align}
The purpose of this paper is to give a $q$-analogue of \eqref{eq:new1} and \eqref{eq:new2} as follows:
\begin{align}
\sum_{k=0}^{\lfloor n/2\rfloor}{m+k\brack k}_{q^2}{m+1\brack n-2k}_{q} q^{n-2k\choose 2}
&={m+n\brack n}_{q}, \label{eq:new3}\\
\sum_{k=0}^{\lfloor n/4\rfloor}{m+k\brack k}_{q^4}{m+1\brack n-4k}_{q} q^{n-4k\choose 2}
&=\sum_{k=0}^{\lfloor n/2\rfloor}(-1)^k{m+k\brack k}_{q^2}{m+n-2k\brack n-2k}_{q}, \label{eq:new4}
\end{align}
where  the {\it $q$-binomial coefficient} ${x\brack k}_{q}$ is defined by
$$
{x\brack k}_{q}=
\begin{cases}
\displaystyle\prod_{i=1}^k\frac{1-q^{x-i+1}}{1-q^i}, &\text{if $k\geqslant 0$},\\[10pt]
0,&\text{if $k<0$.}
\end{cases}
$$
We shall give two proofs of \eqref{eq:new3} and \eqref{eq:new4}. One is combinatorial and the other algebraic.

\section{Bijective proof of \eqref{eq:new3}}
Recall that a {\it partition} $\lambda$ is defined as a finite sequence of
nonnegative integers $(\lambda_1,\lambda_2,\break\ldots,\lambda_r)$
in decreasing order $\lambda_1\geqslant \lambda_2 \geqslant \cdots \geqslant \lambda_r.$
A nonzero $\lambda_i$ is called a {\it part} of $\lambda$. The number of parts of $\lambda$, denoted by $\ell(\lambda)$,
is called the {\it length} of $\lambda$. Write $|\lambda|=\sum_{i=1}^{m}\lambda_i$, called the {\it weight} of $\lambda$.
The sets of all partitions and partitions into distinct parts are denoted by
$\mathscr P$ and $\mathscr D$ respectively. For two partitions $\lambda$ and $\mu$, let $\lambda\cup \mu$
be the partition obtained by putting all parts of $\lambda$ and $\mu$
together in decreasing order.

It is well known that (see, for example, \cite[Theorem 3.1]{Andrews})
\begin{align*}
\sum_{\substack{\lambda_1\leqslant m+1\\ \ell(\lambda)=n}}q^{|\lambda|}
&=q^n{m+n\brack n}_{q} ,  \\
\sum_{\substack{\lambda\in\mathscr D\\\lambda_1\leqslant m+1\\ \ell(\lambda)=n}}q^{|\lambda|}
&={m+1\brack n}_{q} q^{n+1\choose 2}.
\end{align*}
Therefore,
\begin{align*}
\sum_{\substack{\mu\in\mathscr D\\ \lambda_1,\mu_1\leqslant m+1\\ 2\ell(\lambda)+\ell(\mu)=n}}q^{2|\lambda|+|\mu|}
=q^n\sum_{k=0}^{\lfloor n/2\rfloor}{m+k\brack k}_{q^2}{m+1\brack n-2k}_{q} q^{n-2k\choose 2},
\end{align*}
where $k=\ell(\lambda)$. Let
\begin{align*}
\mathscr A&=\{\lambda\in\mathscr P\colon
\lambda_1\leqslant m+1\ \text{and}\ \ell(\lambda)= n\},\\
\mathscr B&=\{(\lambda,\mu)\in\mathscr P\times\mathscr D\colon
\lambda_1,\mu_1\leqslant m+1\ \text{and}\ 2\ell(\lambda)+\ell(\mu)= n\}.
\end{align*}
We shall construct a weight-preserving bijection $\phi$ from $\mathscr A$ to $\mathscr B$.
For any $\lambda\in \mathscr A$, we associate it with a pair $(\overline{\lambda},\mu)$ as follows: If $\lambda_i$ appears $r$ times in $\lambda$,
then we let $\lambda_i$ appear $\lfloor r/2\rfloor$ times in $\overline{\lambda}$ and $r-2\lfloor r/2\rfloor$ times in $\mu$.
For example, if $\lambda=(7,5,5,4,4,4,4,2,2,2,1)$, then
$\overline{\lambda}=(5,4,4,2)$ and $\mu=(7,2,1)$.
Clearly, $(\overline{\lambda},\mu)\in \mathscr B$ and $|\lambda|=2|\overline{\lambda}|+|\mu|$.
It is easy to see that $\phi:\lambda\mapsto(\overline{\lambda},\mu)$ is a bijection. This proves that
$$
\sum_{\lambda\in\mathscr A}q^{|\lambda|}
=\sum_{(\lambda,\mu)\in\mathscr B}q^{2|\lambda|+|\mu|}.
$$
Namely, the identity \eqref{eq:new3} holds.

\section{Involutive proof of \eqref{eq:new4}}
It is easy to see that
\begin{align}
q^n\sum_{k=0}^{\lfloor n/2\rfloor}(-1)^k{m+k\brack k}_{q^2}{m+n-2k\brack n-2k}_{q}
&=\sum_{k=0}^{\lfloor n/2\rfloor}(-1)^k
\sum_{\substack{\lambda_1\leqslant m+1\\ \ell(\lambda)=k}} q^{2|\lambda|}
\sum_{\substack{\mu_1\leqslant m+1\\ \ell(\mu)=n-2k}} q^{|\mu|} \nonumber\\
&=\sum_{\substack{\lambda_1,\mu_1\leqslant m+1\\ 2\ell(\lambda)+\ell(\mu)=n}} (-1)^{\ell(\lambda)} q^{2|\lambda|+|\mu|}, \label{eq:n-2k}
\end{align}
and
\begin{align}
q^n\sum_{k=0}^{\lfloor n/4\rfloor}{m+k\brack k}_{q^4}{m+1\brack n-4k}_{q} q^{n-4k\choose 2}
=\sum_{\substack{\mu\in\mathscr D\\ \lambda_1,\mu_1\leqslant m+1\\ 4\ell(\lambda)+\ell(\mu)=n}}  q^{4|\lambda|+|\mu|}.
\label{eq:n-4k}
\end{align}
Let
\begin{align*}
\mathscr U=\{(\lambda,\mu)\in\mathscr P\times\mathscr P\colon
\lambda_1,\mu_1\leqslant m+1\ \text{and}\ 2\ell(\lambda)+\ell(\mu)= n\}, \\
\mathscr V=\{(\lambda,\mu)\in\mathscr U\colon
\text{each $\lambda_i$ appears an even number of times and $\mu\in\mathscr D$}\}.
\end{align*}
We shall construct an involution $\theta$ on the set $\mathscr U\setminus\mathscr V$
with the properties that $\theta$ preserves $2|\lambda|+|\mu|$ and reverses the sign $(-1)^{\ell(\lambda)}$.

For any $(\lambda,\mu)\in\mathscr U\setminus\mathscr V$, notice that either some $\lambda_i$ appears an odd number of times in $\lambda$,
or some $\mu_j$ is repeated in $\mu$, or both are true. Choose the largest such $\lambda_i$ and $\mu_j$ if they exist, denoted by
$\lambda_{i_0}$ and $\mu_{j_0}$ respectively. Define
\[
\theta((\lambda,\mu))
=\begin{cases}
((\lambda\setminus \lambda_{i_0}),\mu\cup(\lambda_{i_0},\lambda_{i_0})),
&\text{if $\lambda_{i_0}\geqslant \mu_{j_0}$ or $\mu\in\mathscr D$},\\
((\lambda\cup \mu_{j_0}),\mu\setminus(\mu_{j_0},\mu_{j_0})),
&\text{if $\lambda_{i_0}< \mu_{j_0}$ or $\lambda_{i_0}$ does not exist}.
\end{cases}
\]
For example, if $\lambda=(5,5,4,4,4,3,3,3,1,1)$ and $\mu=(5,3,2,2,1)$, then
$$
\theta(\lambda,\mu)=((5,5,4,4,3,3,3,1,1),(5,4,4,3,2,2,1)).
$$
It is easy to see that $\theta$ is an involution on $\mathscr U\setminus\mathscr V$ with the desired properties.
This proves that
\begin{align}
\sum_{(\lambda,\mu)\in\mathscr U} (-1)^{\ell(\lambda)} q^{2|\lambda|+|\mu|}
&=\sum_{(\lambda,\mu)\in\mathscr V} (-1)^{\ell(\lambda)} q^{2|\lambda|+|\mu|}  \nonumber\\
&=\sum_{\substack{\mu\in\mathscr D\\ \tau_1,\mu_1\leqslant m+1\\ 4\ell(\tau)+\ell(\mu)=n}}  q^{4|\tau|+|\mu|}, \label{eq:bijuv}
\end{align}
where $\lambda=\tau\cup\tau$. Combining \eqref{eq:n-2k}--\eqref{eq:bijuv}, we complete the proof of \eqref{eq:new4}.

\section{Generating function proof of \eqref{eq:new3} and \eqref{eq:new4}}
Recall that the {\it $q$-shifted factorial} is defined by
$$(a;q)_0=1,\quad (a;q)_n=\prod_{k=0}^{n-1}(1-aq^{k}),\ n=1,2,\ldots.$$
Then we have
\begin{align}
\frac{1}{(z^2;q^2)_{m+1}}(-z;q)_{m+1}&=\frac{1}{(z;q)_{m+1}}, \label{eq:qbione}\\
\frac{1}{(z^4;q^4)_{m+1}}(-z;q)_{m+1}&=\frac{1}{(z;q)_{m+1}}\frac{1}{(-z^2;q^2)_{m+1}}.  \label{eq:qbitwo}
\end{align}
By the $q$-binomial theorem (see, for example, \cite[Theorem 3.3]{Andrews}), we may expand
\eqref{eq:qbione} and \eqref{eq:qbitwo} respectively as follows:
\begin{align}
&\left(\sum_{k=0}^\infty {m+k\brack k}_{q^2} z^{2k}\right)
\left(\sum_{k=0}^{m+1} {m+1\brack k}_{q} q^{k\choose 2} z^k\right)
=\sum_{k=0}^\infty {m+k\brack k}_{q} z^k, \label{eq:qbie1}\\
&\left(\sum_{k=0}^\infty {m+k\brack k}_{q^4} z^{4k}\right)
\left(\sum_{k=0}^{m+1} {m+1\brack k}_{q} q^{k\choose 2} z^k\right)  \nonumber\\
&=\left(\sum_{k=0}^\infty {m+k\brack k}_{q} z^k\right)
\left(\sum_{k=0}^\infty {m+k\brack k}_{q^2} (-1)^k z^{2k}\right).  \label{eq:qbie2}
\end{align}
Comparing the coefficients of $z^n$ in both sides of \eqref{eq:qbie1} and \eqref{eq:qbie2},
we obtain \eqref{eq:new3} and \eqref{eq:new4} respectively.

Finally, we give the following special cases of \eqref{eq:new3}:
\begin{align}
\sum_{k=0}^{\lfloor n/2\rfloor}{n+k\brack k}_{q^2}{n+1\brack 2k+1}_{q} q^{n-2k\choose 2}
&={2n\brack n}_{q}, \label{eq:spe1} \\
\sum_{k=0}^{\lfloor n/2\rfloor}{n+k\brack k+1}_{q^2}{n\brack 2k+1}_{q} q^{n-2k-1\choose 2}
&={2n\brack n-1}_{q}.  \label{eq:spe2}
\end{align}
When $q= 1$, the identities \eqref{eq:spe1} and \eqref{eq:spe2} reduce to \eqref{eq:s1} and \eqref{eq:s2} respectively.

\vskip 5mm
\noindent{\bf Acknowledgments.} This work was partially
supported by the Fundamental Research Funds for the Central Universities,  Shanghai Rising-Star Program (\#09QA1401700),
Shanghai Leading Academic Discipline Project (\#B407), and the National Science Foundation of China (\#10801054).

\end{document}